\documentclass[11pt]{article}
\usepackage{fullpage}

\pdfoutput=1
\usepackage{graphicx, rotating}
\usepackage{latexsym,amsmath,amssymb,amsthm,epic,eepic,multirow}
\usepackage{natbib}
\usepackage{verbatim} 
\usepackage{url}
\usepackage{setspace}
\usepackage{graphicx}
\usepackage{multirow}
\usepackage{subfigure}
\usepackage{amsfonts}
\usepackage{natbib}
\usepackage{algorithm}
\usepackage{algorithmic}
\usepackage[pdftex]{color}
\usepackage{colortbl}
\definecolor{first}{rgb}{0.65,0.65,0.65}
\definecolor{second}{rgb}{0.9,0.9,0.9}
\definecolor{brgreen}{rgb}{0.0, 0.26, 0.15}
\definecolor{cadmiumgreen}{rgb}{0.0, 0.42, 0.24}

\newcommand{\PP}{\mathbb{P}}
\newcommand{\EE}{\mathbb{E}}

\newcommand{\eps}{\varepsilon}

\newcommand{\bx}{\mathbf{X}}
\newcommand{\by}{\mathbf{Y}}

\newcommand{\rank}[1]{\mbox{rank}}

\newcommand{\NN}{\mathbb{N}}
\newcommand{\si}{\mbox{si}}
\newcommand{\ch}{\mbox{ch}}
\newcommand{\pa}{\mbox{pa}}

\newcommand{\aggr}{\mbox{aggr}}

\newcommand{\anc}{\mbox{anc}}
\newcommand{\of}{\mbox{of}}

\newtheorem{theo}{Theorem}
\newtheorem{prop}{Proposition}

\usepackage{xr}

\begin{document}
\title{A sequential rejection testing method for\\
high-dimensional regression with correlated variables} 

\author{}
\author{Jacopo Mandozzi and Peter B\"uhlmann\\
Seminar for Statistics, ETH Z\"urich}

\maketitle

\begin{abstract}
We propose a general, modular method for significance testing of groups
(or clusters) of variables in a high-dimensional linear model. In presence
of high correlations among the covariables, due to serious problems of
identifiability, it is indispensable to focus on detecting
groups of variables rather than singletons.   
We propose an inference method which allows to build in hierarchical
structures. It relies on repeated sample
splitting and sequential rejection, and we prove that it asymptotically
controls the familywise error rate. It can be implemented on any
collection of clusters and leads to improved power in comparison to more
standard non-sequential rejection methods. 
We complete the theoretical analysis with empirical results for simulated
and real data. 
 
\end{abstract}

{\bf Keywords and phrases:} Familywise error rate; Hierarchical
clustering; High-dimensional
variable selection; Inheritance procedure; Lasso;
Linear model; Minimal true detection; Multiple testing; Sample
splitting; Sequential rejection principle; Singleton true detection.

\section{Introduction}
Error control of false selection or false positive statements based on
p-values is a primary goal of
statistical inference and an established, broadly used tool in many areas
of science. It relies on standard statistical hypothesis testing and
procedures which give provable guarantees in presence of multiple,
potentially very large scale multiple testing
\citep{westfall1993resampling,dudoit2007multiple,efron2010large}. While
being standard in 
the classical low-dimensional setup, statistical significance testing in
the more challenging high-dimensional setting where the number of variables
$p$ might be much larger than the sample size $n$ has only received
attention recently.  

We consider here a linear regression model 
\begin{eqnarray}\label{mod1}
\by = \bx \beta^0 + \eps,\ \eps \sim {\cal N}_n(0,\sigma^2 I),
\end{eqnarray}
with $n \times p$ design matrix $\bx$, $p \times 1$ regression vector
$\beta^0$ and $n \times 1$ response $\by$. We allow for high-dimensional
scenarios where $p \gg n$. We assume that the regression coefficient vector
is sparse with many coefficients of $\beta^0$ being equal to zero, that is,
the active set of variables $$S_0 = \{ j;
\beta^0_j \neq 0 \}$$ is assumed to be a small subset of $\{1,\ldots ,p\}$
corresponding to all variables. 

A few methods for assigning $p$-values and
constructing confidence intervals for individual parameters $\beta^0_j\
(j=1,\ldots ,p)$ have been suggested
\citep{WR08,memepb09,bue12,zhazha13,geeretal13,loc14,javmo13}, and some of
them have been compared against each other in various settings
\citep{bumeka13, DeBuMeMe14}. The inferential statements can easily
be adjusted for multiplicity, thanks to the methodology 
and theory in multiple testing \citep[cf.]{dudoit2007multiple}. However, and
important for practical 
applications, some major issues in presence of highly correlated variables
still need further attention: typically, when $p \gg n$, none or only a few
of the individual $\beta^0_j$'s turn out to be significant which is a
consequence of their near non-identifiability (even when some theoretical
conditions on well-posedness on the design matrix $\bx$
\citep[cf.]{pbvdg11} hold). However, a \emph{group} of (correlated) variables is
often much better identifiable, but one can then not determine anymore 
the relevant variables within such a group \citep{buru12,meins13,MaBu14}.    
 
Thus, our main goal is testing of significance of groups of parameters: 
for a group or cluster $C \subseteq \{1,\ldots ,p\}$ we consider the
following null- and alternative hypothesis, respectively: 
\begin{eqnarray*}
H_{0,C}:\ \beta^0_j = 0\ \mbox{for all}\ j \in C,\ \ \ H_{A,C}:\ \beta^0_j
\neq 0\ \mbox{for at least one}\ j \in C.
\end{eqnarray*}

Given a collection $\mathcal{C}$ of clusters, we propose a
general method for obtaining a 
collection $\mathcal{R} \subset 
\mathcal{C}$ of rejected clusters such that familywise error rate (FWER) is
strongly controlled. That is, for a given nominal level $\alpha \in (0,1)$:
$$\PP[\mathcal{R} \subseteq \mathcal{F}]  \geq 1-\alpha,$$
where $\mathcal{F} = \{ C \in \mathcal{C} \mbox{ s.t. }  H_{0,C}
\mbox{ is false}\}$ i.e., ${\mathcal F}$ is the collection of false null
hypotheses.
Our new method has the following main features:
\begin{itemize}
\item It can be implemented on any collection of clusters $\mathcal{C}$.
\item It is modular in the sense that it requires four basic building
  blocks that have to satisfy certain assumptions.
\item Its modular conception allows for a better insight of the procedure's
  power and improvements thereof.
\end{itemize}

We are particularly interested to use the procedure for
\emph{hierarchically ordered} clusters of (correlated) variables. Such a
hierarchical structure can be obtained from the output of a hierarchical
clustering algorithm: since it operates on the design matrix $\bx$ only and
does not involve the responses $\by$, the inference for $\beta^0$ remains
correct (for fixed design or by conditioning on $\bx$). With such a
hierarchical cluster 
tree, our inference method (Sections \ref{sec:procedure} and
\ref{sec:hierarchical}) first tests the cluster ${\cal C} = \{1,\ldots ,p\}$
containing all the 
variables (the top node in the tree): if the
corresponding null-hypothesis is rejected, we test some refined clusters, and
we proceed down the cluster tree, in a \emph{sequential manner}, until a
cluster is not significant anymore. Figures \ref{figure:dendro1} and
\ref{figure:dendro2} in Section \ref{sec:empirical} provide some graphical
illustrations. This procedure has 
the remarkable property that the 
resolution level of the significant clusters is automatically controlled by the
sequential testing method: if the signal is strong (e.g. large absolute
values of components of $\beta^0$) and the variables are not too highly
correlated, one can detect small clusters or even single variables and
vice-versa, if the signal isn't very strong or the variables are highly
correlated, only larger groups can be detected as significant. 

\paragraph{Relation to other work.}
Our proposed method is based on the multi sample splitting method from
\citet{memepb09} and the 
sequential rejection principle of \citet{GoSo2010}. It is a generalization
and power improvement over the multi sample splitting technique for
inference of single variables \citep{memepb09} and for hierarchically
ordered clusters of variables \citep{MaBu14}. The improvement in power is
strict, and in analogy to the gain of power of Holm's procedure \citep{Holm1979} over the Bonferroni adjustment. Thus, even if
the increased power might be only small for some datasets, one cannot do worse
with the new procedure. The only price to pay is a slightly more
complicated algorithm: we provide an implementation in the \textrm{R}-package
\texttt{hdi}. 

\paragraph{Outline of the paper.} 
In Section \ref{sec:buildingblocks} we describe the four basic building blocks of the method and the assumptions that are
sufficient to establish in Section \ref{sec:FWER} its strong FWER control.
In Sections \ref{sec:single} and \ref{sec:hierarchical}, respectively, we
focus on the inference of two specific kinds of  
cluster collections: singletons and hierarchically ordered clusters. In
Section \ref{sec:Shaffer} we show how logical relationships can be used to
improve the power. Finally, we provide in Section \ref{sec:empirical} a
comparison based on empirical results for error control and power, with a
focus on minimal true detections, and we apply the new method to a real
dataset. 

\section{A construction based on four building
  blocks}\label{sec:buildingblocks}  
Our method is based on four basic building blocks that satisfy certain
assumptions. 

One main ingredient is multi sample splitting. For $b=1,\ldots, B$ where
$B$ is the number of repeated sample splitting, the original data of sample
size $n$ is split into two disjoint groups, $N_{in}^{(b)}$ and
$N_{out}^{(b)}$, i.e., a partition 
$$\{1,\dots,n \} = N_{in}^{(b)} \cup N_{out}^{(b)}$$ is randomly chosen. The groups
are chosen of equal size if $n$ is even or satisfy
$|N_{out}^{(b)}|=|N_{in}^{(b)}|+1$ if $n$ is odd. 

The idea is to use data from $N_{in}^{(b)}$ to select a few variables and
the other data from $N_{out}^{(b)}$ to perform the statistical hypothesis
testing in the low-dimensional submodel with the selected variables from
$N_{in}^{(b)}$. The details are described next. 

\subsection{Screening of variables}\label{sec:screening}
We consider variable screening where an estimator
$\hat{S}^{(b)} \subseteq \{1,\ldots ,p\}$, based on data corresponding to
$N_{in}^{(b)}$, is aiming at including all active
variables $S_0$. A prime example is the Lasso \citep{T96}, while a detailed
empirical comparison of five popular screening procedures can be found in \citep{bueman12}. 
%
Assume that the screening procedure satisfies the following
properties for any sample split $b$:
\begin{eqnarray*}
&&\mbox{(A1) \it{Sparsity property}:}\; |\hat{S}^{(b)}| < n/2.\\
&&\mbox{(A2) \it{$\delta$-Screening property}:}\; \PP[\hat{S}^{(b)} \supseteq
S_0] \geq 1-\delta,\ \mbox{where}\ 0 < \delta < 1.
\end{eqnarray*}
The \textit{sparsity property} in (A1) implies that for each sample split
$b$ it holds that $|\hat{S}^{(b)}|<|N_{out}^{(b)}|$, a condition which is
necessary for applying classical tests as described in Section
\ref{sec:testing} below. The \textit{$\delta$-screening
  property} in (A2) ensures that all the 
relevant variables are retained with high probability (where $\delta > 0$ is
typically small).

We indicate in Section \ref{sec:justific} that under some assumptions, the
Lasso satisfies (A1) and (A2). 

\subsection{Testing and p-values}\label{sec:testing}

The idea is to perform a classical statistical test on the other half
sample from $N_{out}^{(b)}$ in a low-dimensional problem with variables from
$\hat{S}^{(b)}$ only.  

For each sample split $b$, based on the second half of the sample
corresponding to $N_{out}^{(b)}$, consider a testing
procedure, e.g. the classical partial F-test (see also Section
\ref{sec:justific}), that provides correct 
p-values $p^{C,(b)}$ for the null hypothesis $H_{0,C \cap \hat{S}^{(b)}}$ for
each screened set $\hat{S}^{(b)}$, in the sense that for each nominal 
level $\alpha \in (0,1)$ 
\begin{eqnarray*}
\mbox{(A3) \it{Correct testing property}:} \mbox{ Under the null hypothesis
} H_{0,C \cap \hat{S}^{(b)}}\, \mbox{it holds } \PP[p^{C,(b)}\leq \alpha] \leq \alpha.
\end{eqnarray*}
We note that the probability is with respect to the data generating random
variables corresponding to the second half $N_{out}^{(b)}$, and the
null-hypothesis is fixed with respect to $N_{out}^{(b)}$. Due to the
screening property (A2), when $\delta \to 0$, the null-hypothesis $H_{0,C \cap
  \hat{S}^(b)}$ approximates the unconditional hypothesis $H_{0,C}$ which
we aim to test for. If $C \cap \hat{S}^{(b)} = \emptyset$ define
$p^{C,(b)}=1$. This provides a 
(correct) p-value $p^{C,(b)}$ for each cluster $C \in \mathcal{C}$ and each
sample split $b \in \{1 \dots B\}$.

\subsection{Multiplicity adjustment}\label{sec:multadj}
Consider for each sample split $b$ and each cluster 
$C \in \mathcal{C}$ a multiplicity adjustment procedure 
$m^{(b)}_C:2^{\mathcal{C}} \rightarrow [1,\infty]$ that for each collection
$\mathcal{R}$ of rejected 
clusters provides a multiplicity adjustment
$m^{(b)}_C(\mathcal{R}) \geq 1$
and satisfies the following properties:
\begin{eqnarray*}
&&\mbox{(A4) \it{Monotonicity property}:}\; \mbox{If } \mathcal{R} \subseteq \mathcal{S} \mbox{ then } m^{(b)}_C(\mathcal{R}) \geq m^{(b)}_C(\mathcal{S}).\\
&&\mbox{(A5) \it{Single-step property}:}\; \sum_{C \in\, \mathcal{C} \setminus
  \mathcal{R}} \frac{1\{C \cap \hat{S}^{(b)} \neq \emptyset\}}{m^{(b)}_C(\mathcal{R})} \leq 1, 
\end{eqnarray*}
where we define $1/\infty=0$. Such a family of multiplicity adjustments for
$b=1,\ldots, B$ are often naturally induced from a global multiplicity
adjustment procedure $m_C$. 


\subsection{Aggregation of p-values}\label{sec:aggr}

Consider a collection of screened sets of variables $\hat{S}^{(b)}$, a
cluster $C \in \mathcal{C}$, a collection of p-values $p^{C,(b)}$ for the
null-hypothesis $H_{0,C \cap S^{(b)}}$ (which approximates $H_{0,C}$, see
comment after (A3)) and a collection of multiplicity adjustments
$m_C^{(b)}\geq1$ (we drop here the dependence on ${\cal R}$). 

The goal is to aggregate the p-values $p^{C,(1)},\ldots , p^{C,(B)}$ to a
single p-value which is adjusted for multiplicity. 
An  aggregation procedure is a monotone increasing function $\aggr: [0,1]^B 
\rightarrow [0,1]$. Assume it satisfies the following property:
\begin{eqnarray*}
\mbox{(A6) \it{Aggregation property}:}&& \mbox{If } \PP[p^{C,(b)} \leq \alpha]
\leq \alpha, \ \forall \alpha \in [0,1], \mbox{ then}\\
&&\PP[\aggr(p^{C,(1)}m_C^{(1)},\dots,p^{C,(B)}m_C^{(B)}) \leq \alpha]\\
&&\leq \frac{\alpha}{B} \sum_{b=1}^B \frac{1\{C \cap \hat{S}^{(b)} \neq 
  \emptyset\}}{m_C^{(b)}}, \ \forall \alpha \in [0,1]. 
\end{eqnarray*}

\subsection{The procedure}\label{sec:procedure}

Our procedure is based on the four building blocks above. First, we proceed
with screening of the variables based on the first half sample from
$N_{in}^{(b)}$ (Section \ref{sec:screening}), e.g., in Section
\ref{sec:impscen} we use the Lasso with regularization parameter
chosen by 10-fold cross-validation (see also Section \ref{sec:justific}). 
Then, we construct the p-values based on the second half sample from
$N_{out}^{(b)}$ by using 
the partial F-test (Section \ref{sec:testing} and see also Section
\ref{sec:justific}). This leads to a (correct) p-value $p^{C,(b)}$ for each 
cluster $C \in \mathcal{C}$ and each sample split $b \in \{1 \dots
B\}$. 

The multiplicity adjustment is done sequentially (Section
\ref{sec:multadj}). Based on a chosen significance level $\alpha \in (0,1)$
and for a collection of currently rejected sets $\mathcal{R}$, define the
successor of $\mathcal{R}$ as 
$$\mathcal{N}(\mathcal{R})=\{ C \in \mathcal{C} \setminus \mathcal{R}
\mbox{ s.t. }
\aggr(p^{C,(1)}m^{(1)}_C(\mathcal{R}),\dots,p^{C,(B)}m^{(B)}_C(\mathcal{R}))
\leq \alpha \}$$  
Start from ``no rejections'' $\mathcal{R}_0 = \emptyset$, define
$\mathcal{R}_{i+1}=\mathcal{R}_i \cup \mathcal{N}(\mathcal{R}_i)$ and 
$\mathcal{R}_\infty = \lim_{i \to \infty} \mathcal{R}_i$ (although
  ${\cal R}_{\infty}$ is never constructed due to finite-ness of all possible
  subset of the variables). 
Concrete choices of $m^{(1)}_C(\mathcal{R}),\ldots ,m^{(B)}_C(\mathcal{R})$
are discussed in Section \ref{sec:multiplicity}.

Finally, we aggregate the p-values as indicated in Section
\ref{sec:aggr}. Concrete aggregation methods are described in Proposition
\ref{theo:quantiles} in Section \ref{sec:justific}.

\section{Familywise error control}\label{sec:FWER}

We show here that the method from Section \ref{sec:procedure} (strongly)
controls the FWER at each step $i=0,1,2,\ldots $ 
\begin{theo}\label{theo:FWER}
Assume that (A1)-(A6) hold. Then for any $i \in \NN_0 \cup \infty$
$$\PP[\mathcal{R}_i \subseteq \mathcal{F}] \geq (1-\delta)^B-\alpha,$$
  where $\mathcal{F} = \{ C \in \mathcal{C} \mbox{ s.t. } C \cap S_0 \neq
\emptyset\}$ is the collection of false null hypotheses.
\end{theo}
A proof is given in the Appendix. 

\subsection{Screening, testing and aggregation: their properties}\label{sec:justific}

We discuss here some choices for screening, testing and aggregation which
we use in the implementation in the \textrm{R}-package \texttt{hdi}. The
issue of sequential multiplicity adjustment is treated separately in Section
\ref{sec:multiplicity}. 

For variable screening, we use the Lasso with regularization parameter
chosen by 10-fold cross-validation. Theoretical justification of the
sparsity and screening property (A1) and
(A2) can be derived by assuming a compatibility or restricted eigenvalue
condition on the fixed design 
matrix $\bx$ and a beta-min assumption requiring that $\min_{j \in
  S_0}|\beta^0_j| \gg \sqrt{|S_0| \log(p)/n}$ is sufficiently large: we
refer to \citet[Ch. 2.7 and Ch. 6]{pbvdg11} for the details. 

For construction of the p-values (in the low-dimensional setting, due to
variable screening in the first half of the sample) we use the partial
F-test. Then, assuming fixed design $\bx$ and Gaussian errors, condition
(A3) holds. 

For aggregation of the p-values, ensuring that (A6) holds, we have the
following result for two slightly different methods.
 
\begin{prop}\label{theo:quantiles}
Denote by $q_{\gamma}(u)$ the empirical $\gamma$-quantile of the values
occurring in the components of a vector $u$. The monotone increasing
functions $[0,1]^B \rightarrow [0,1]$
\begin{eqnarray*}
\big( \tilde{p}^{(1)},\dots,\tilde{p}^{(B)} \big) &\longmapsto& Q(\gamma)=\min \big\{ \,1~,~q_\gamma
\big( \tilde{p}^{(1)}/\gamma,\dots,\tilde{p}^{(B)}/\gamma \big) \big\}\\
\big( \tilde{p}^{(1)},\dots,\tilde{p}^{(B)} \big) &\longmapsto& \min \big\{ \,1~,~(1-\log\gamma_{\min}) \inf_{\gamma \in (\gamma_{\min},1)} Q(\gamma) \big\}
\end{eqnarray*}
satisfy the aggregation property (A6) for any $\gamma,\gamma_{min} \in (0,1)$.  
\end{prop}
A proof, which was basically given in \citet{memepb09}, can be found in
the Appendix. 

\section{Some concrete methods for multiplicity adjustment}\label{sec:multiplicity}

We discuss here the issue of multiplicity adjustment, and justify
assumption (A4) and (A5) for different inference procedures. 

\subsection{Inference of single variables}
\label{sec:single}
This first example is paradigmatic for the advantages of the modular
approach: a simple improvement of the multiplicity adjustment procedure
allows for a better power, basically in the same way as in a
low-dimensional setting in \citep{GoSo2010}. 

Concretely, we consider the problem of inferring single variables, i.e.,
the collection of clusters $\mathcal{C} = \{ \{i\}; i=1,\dots,p\}$. The 
method proposed in 
\citet{memepb09} corresponds to the method of Theorem \ref{theo:FWER} with
the aggregation procedures of Proposition \ref{theo:quantiles} and the
following Bonferroni-based \citep{Bonf36, Dunn61} multiplicity adjustment
procedure:
\begin{equation}\label{eq:originalsingle}
m^{(b)}_{\{i\}}(\mathcal{R}) = |\hat{S}^{(b)}|.
\end{equation}
As the multiplicity adjustments are independent from the (previously)
rejected collection of sets, the monotonicity property (A4) is trivially
satisfied, while the single-step property (A5) follows from
\begin{eqnarray*}
\sum_{C \in\, \mathcal{C} \setminus
  \mathcal{R}} \frac{1\{C \cap \hat{S}^{(b)} \neq \emptyset\}}{m^{(b)}_C(\mathcal{R})} = \sum_{\{i\} \in\, \mathcal{C} 
\setminus \mathcal{R}} \frac{1\{\{i\} \cap \hat{S}^{(b)} \neq
\emptyset\}}{|\hat{S}^{(b)}|}\leq 1. 
\end{eqnarray*}
The power of the method can be improved taking instead of
(\ref{eq:originalsingle}) the following Bonferroni-Holm-based
\citep{Holm1979} multiplicity adjustment procedure:
\begin{equation}\label{eq:Holmsingle}
m^{(b)}_{\{i\}}(\mathcal{R}) = |\hat{S}^{(b)} \cap (\mathcal{C} \setminus
\mathcal{R})|=|\{j \in \hat{S}^{(b)} \mbox{ s.t. } \{j\} \notin \mathcal{R}\}|.
\end{equation}
The monotonicity property (A4) is still satisfied since $|\hat{S}^{(b)} \cap
(\mathcal{C} \setminus \mathcal{R})| \geq |\hat{S}^{(b)} \cap (\mathcal{C}
\setminus \mathcal{S})|$ for $\mathcal{R} \subseteq \mathcal{S}$, whereas  
\begin{eqnarray*}
\sum_{C \in\, \mathcal{C} \setminus
  \mathcal{R}} \frac{1\{C \cap \hat{S}^{(b)} \neq \emptyset\}}{m^{(b)}_C(\mathcal{R})} = \sum_{\{i\} \in\, \mathcal{C} 
\setminus \mathcal{R}} \frac{1\{\{i\} \cap \hat{S}^{(b)} \neq
\emptyset\}}{|\hat{S}^{(b)} \cap (\mathcal{C} \setminus \mathcal{R})|} = 1 
\end{eqnarray*}
proves the single step property (A5). 

\subsection{Inference of hierarchically ordered clusters of variables}\label{sec:hierarchical}
When dealing with the challenge of inferring hierarchically ordered
clusters of variables, e.g. from the tree-structured output of a
hierarchical clustering 
algorithm, one considers a collection of clusters $\mathcal{C} =
\{C_i\}_i$ where for any  two clusters $C_i,C_{i'} \in \mathcal{C}$, either
one cluster is a subset of the other, or they have an empty intersection. The
method proposed in \citet[Section 2]{MaBu14}, which is based on the procedure of
\citet{Meins08}, corresponds to the one as in Theorem \ref{theo:FWER} with
the aggregation methods of Proposition \ref{theo:quantiles} and the
following multiplicity adjustment:
\begin{equation}\label{eq:originalhier}
m^{(b)}_{C}(\mathcal{R}) = 
\left\{ 
\begin{array}{ll}
\infty, & \mbox{if } \anc(C) \not\subseteq \mathcal{R}\\
\frac{|\hat{S}^{(b)}|}{|\hat{S}^{(b)} \cap C|}, & \mbox{if } \anc(C)
\subseteq \mathcal{R} \mbox{ and }\hat{S}^{(b)} \cap C \neq \emptyset\\
1, & \mbox{otherwise.}
\end{array}\right.
\end{equation}
Here, $\anc(C)$ denotes the ancestors in a hierarchically ordered cluster
tree. To check the monotonicity property (A4), consider $\mathcal{R} \subseteq
\mathcal{S}$. For $C \in \mathcal{C}$ with $\anc(C)
\subseteq \mathcal{R}$ it holds $\anc(C)
\subseteq \mathcal{S}$ and hence $m^{(b)}_{C}(\mathcal{R}) =
m^{(b)}_{C}(\mathcal{S})$, while for $C \in \mathcal{C}$ with $\anc(C)
\not\subseteq \mathcal{R}$ one has $m^{(b)}_{C}(\mathcal{R}) = \infty \geq
m^{(b)}_{C}(\mathcal{S})$. The single step property (A5) follows from
\begin{eqnarray*}
\sum_{C \in\, \mathcal{C} \setminus
  \mathcal{R}} \frac{1\{C \cap \hat{S}^{(b)} \neq
  \emptyset\}}{m^{(b)}_C(\mathcal{R})} &=& \frac{1}{|\hat{S}^{(b)}|} \sum_{C \in\, \mathcal{C} \setminus
  \mathcal{R} \mbox{ s.t. } \anc(C) \,\subseteq\, \mathcal{R}}
|\hat{S}^{(b)} \cap C| \leq 1,
\end{eqnarray*}
where in the inequality we have used the fact that for two sets in the sum
above, one cannot be a subset of the other and hence, by
definition of the hierarchy $\mathcal{C}$, they are
disjoint.

\subsubsection{The inheritance procedure in the high-dimensional setting}\label{sec:inher}
In \citet[Section 6.3]{GoSo2010} and \citet{GoFi2012}, the authors propose
various possibilities on how the sequential rejection principle can be used to
improve the power of the hierarchical procedure in 
\citet{Meins08}. 
We
consider here the most powerful one, the inheritance procedure of
\citet{GoFi2012} which we extend to the high-dimensional setting with
hierarchical cluster trees. In order to
do that, we apply the method of Theorem \ref{theo:FWER} with the aggregation
procedures of Proposition \ref{theo:quantiles} and the following multiplicity
adjustment:
\begin{equation}\label{eq:GFhier}
m^{(b)}_{C}(\mathcal{R}) = 
\left\{ 
\begin{array}{ll}
\infty, & \mbox{if } \anc(C) \not\subseteq \mathcal{R}\\
1, & \mbox{if } \anc(C)
\subseteq \mathcal{R} \mbox{ and }\hat{S}^{(b)} \cap C = \emptyset\\
\frac{|\hat{S}^{(b)}|}{|\hat{S}^{(b)} \cap C|} \prod_{D \in \anc(C)} n_D^{(b)}({\mathcal{R}}) , & \mbox{otherwise,}
\end{array}\right.
\end{equation}
where 
$$n_D^{(b)}(\mathcal{R}) = \frac{1}{|\hat{S}^{(b)} \cap D|} \sum_{E \in \,\ch(D)
  \setminus \mathcal{E}(\mathcal{R})} |\hat{S}^{(b)} \cap E|$$ and 
$$\mathcal{E}(\mathcal{R}) = \{C \in \mathcal{C} \mbox{ s.t. } \of(C)
\subseteq{R} \}$$
the set of extinct branches, i.e., the set of hypotheses which have been
rejected together with all their offsprings denoted by $\of(C)$ (as before,
$\anc(C)$ denotes the ancestors of cluster $C$). Note that since
$n_D^{(b)}({\mathcal{R}}) \leq 1$ this procedure leads to an uniform improvement
over the method of the previous section.

The monotonicity property (A4) follows from the same considerations as
above and
$$ \mathcal{R} \subseteq \mathcal{S} \Longrightarrow \mathcal{E}(\mathcal{R})
\subseteq \mathcal{E}(\mathcal{S}) \Longrightarrow
n_D^{(b)}({\mathcal{R}}) \geq n_D^{(b)}({\mathcal{S}}).$$
To check that the single step property (A5) holds note that
\begin{eqnarray*}
&&\sum_{C \in\, \mathcal{C} \setminus
  \mathcal{R}} \frac{1\{C \cap \hat{S}^{(b)} \neq
  \emptyset\}}{m^{(b)}_C(\mathcal{R})}\\ 
&=& \sum_{C \in\, \mathcal{C} \setminus
  \mathcal{R} \mbox{ s.t. } \anc(C) \,\subseteq\, \mathcal{R}}
\frac{|\hat{S}^{(b)} \cap C|}{|\hat{S}^{(b)}|} \prod_{D \in \anc(C)}
\frac{|\hat{S}^{(b)} \cap D|}{\sum_{E \in \ch(D) \setminus
    \mathcal{E}(\mathcal{R})}|\hat{S}^{(b)} \cap E|}\\
&=& \sum_{C \in\, \mathcal{C} \setminus
  \mathcal{R}} \frac{\alpha^{(b)}_C(\mathcal{R})}{\alpha}
\end{eqnarray*}
where $\alpha^{(b)}_C$ is as in \citet[equation (5)]{GoFi2012} with the
weights $w^{(b)}_C=|\hat{S}^{(b)} \cap C|$; therefore the single step property
follows directly from the considerations in \citet{GoFi2012}.

\subsection{Exploiting logical relationships: Shaffer
  improvements}\label{sec:Shaffer}

Logical relationships between hypothesis can be exploited to improve the
power of the sequential rejection procedure. A first example of such an
improvement for hierarchically 
ordered clusters was given in \citet{Meins08}, while in \citet{GoFi2012}
the improvement is applied to the inheritance procedure. Since those
improvements are based on the considerations of \citet{Shaf86} they
are called ``Shaffer improvements''. For the high-dimensional setting a
possible Shaffer improvement consists of multiplying the multiplicity adjustment
$m_C^{(b)}(\mathcal{R})$ with the Shaffer factor 
\begin{equation}\label{eq:Shaffactor}
s_C^{(b)}(\mathcal{R}) = \max \{ m_C^{(b)}(\mathcal{U}) /
m_C^{(b)}(\mathcal{R})\, \mbox{ s.t. } C \notin \mathcal{U} \supseteq
\mathcal{R},\ \mathcal{U} \mbox{ congruent}\},
\end{equation}
where a set $\mathcal{U} \subseteq \mathcal{C}$ is called congruent if, by
the logical implications, it can be a complete set of false hypothesis
(e.g. for a collection $\mathcal{C}$ of hierarchically ordered hypothesis
$\mathcal{U} \subseteq \mathcal{C}$ is congruent if for each $C \in
\mathcal{U}$ it holds $\anc(C) \subseteq \mathcal{U}$ and at least one offspring
leaf node of $C$ is in $\mathcal{U}$).\\

Note that multiplication with the Shaffer factor never decreases the power of the
method since by the monotonicity property (A4), $s_C^{(b)}(\mathcal{R}) \leq
1$. Moreover $s_C^{(b)}(\mathcal{R}) = 1$ if $\mathcal{R}$ is congruent and
since the collection $\mathcal{F}$ of all false hypothesis is congruent,
the Shaffer improvement doesn't affect the validity of equation
(\ref{eq:step}). Finally, for $\mathcal{R} \subseteq \mathcal{S}$,
\begin{eqnarray*}
m_C^{(b)}(\mathcal{R})s_C^{(b)}(\mathcal{R}) &=& \max \{
m_C^{(b)}(\mathcal{U}) \mbox{ s.t. } C \notin \mathcal{U} \supseteq
\mathcal{R},\ \mathcal{U} \mbox{ congruent}\}\\ 
&\geq& \max \{
m_C^{(b)}(\mathcal{U}) \mbox{ s.t. } C \notin \mathcal{U} \supseteq
\mathcal{S},\ \mathcal{U} \mbox{ congruent}\}\\ 
&=& m_C^{(b)}(\mathcal{S})s_C^{(b)}(\mathcal{S})
\end{eqnarray*}
and hence the Shaffer improvement doesn't affect the validity of equation
(\ref{eq:mono}) neither.\\

We want to apply this Shaffer improvement to the inheritance procedure
described in Section \ref{sec:inher}. Following the same reasoning as in
\citet[Section 6]{GoFi2012}, with the weights $w_C^{(b)}=|\hat{S}^{(b)}
\cap C|$ we get the Shaffer factor
$$s_C^{(b)}(\mathcal{R}) =
\left\{\begin{array}{ll}
\frac{w_C^{(b)}+u_C^{(b)}-v_C^{(b)}}{w_C^{(b)}+u_C^{(b)}}, & \mbox{if } C \notin \mathcal{R}, \si(C) \subseteq \mathcal{L}
\setminus \mathcal{R}\\
1, & \mbox{otherwise,}
\end{array}\right.$$
where $\si(C)=\ch\{\pa(C)\} \setminus \{C\}$ denotes the siblings of $C$,
$\mathcal{L} \subset \mathcal{C}$ denotes the collection of leaf nodes,
$u_C^{(b)}=\sum_{D \in \si(C)}w_D^{(b)}$ and $v_C^{(b)}=\min_{D \in
  \si(C)}w_D^{(b)}$. If $\mathcal{C}$ is a binary tree the Shaffer factor becomes
$$s_C^{(b)}(\mathcal{R}) =
\left\{\begin{array}{ll}
\frac{|\hat{S}^{(b)}\cap C|}{|\hat{S}^{(b)}\cap C|+|\hat{S}^{(b)} \cap
  \,\si(C)|}, & \mbox{if } C \notin \mathcal{R}, \si(C) \subseteq \mathcal{L}
\setminus \mathcal{R}\\
1, & \mbox{otherwise.}
\end{array}\right.$$
Unlike as for the inheritance procedure in (\ref{eq:GFhier}), the
Shaffer factor (\ref{eq:Shaffactor}) for the procedure 
in (\ref{eq:originalhier}) is always 1. Nevertheless, a possibility how to
exploit logical relationships to improve the power of the procedure
(\ref{eq:originalhier}) for binary trees, which provides a Shaffer
improvement very similar to the one above, is illustrated in \citet{MaBu14}.

\section{Empirical results}\label{sec:empirical}

\subsection{Implementation of the methods and considered scenarios}\label{sec:impscen}
In this section we compare the performance of the four methods illustrated in
Sections \ref{sec:single} and \ref{sec:hierarchical} and refined in
Section \ref{sec:Shaffer}, i.e. single variable method with Bonferroni
multiplicity adjustment (\ref{eq:originalsingle}), hierarchical method with Bonferroni-based adjustment
(\ref{eq:originalhier}) along with Shaffer improvement as in
\citet{MaBu14}, single variable method with
Bonferroni-Holm multiplicity adjustment (\ref{eq:Holmsingle}) and
hierarchical method with inheritance procedure (\ref{eq:GFhier}) along with
Shaffer improvement (\ref{eq:Shaffactor}). In the following we refer
  to the first two methods as the ``non-sequential methods'' (strictly
  seen, the hierarchical method with Bonferroni-based adjustment is
  actually sequential, but there previous rejections are not used to
  improve subsequent multiplicity corrections) and the latter
  two methods as the ``sequential methods''.

We consider the same
implementation of the methods and the same scenarios (with exactly
the same sample splits) as in \citet{MaBu14}, although here we use only standard
hierarchical clustering for the hierarchical
methods. Concretely, the following choices have been made for
implementation:
\begin{itemize}
\item construction of the clusters with standard hierarchical clustering
  (using the 
  \texttt{R}-function \texttt{hclust}) with distance between two covariables
  equal to 1 minus the absolute correlation between the covariables, and using
  complete linkage;
\item screening with the Lasso \citep{T96} with regularization parameter
  chosen by 10-fold cross-validation;
\item $B=50$ sample splits (for each scenario exactly the same splits as in
  \citet{MaBu14});
\item for aggregation, the p-values $P_h^C$ in Proposition \ref{theo:quantiles}
  are computed over a grid of  $\gamma$-values between $\gamma_{min}=0.05$
  and $1$ with grid-steps of size $0.025$;
\item nominal significance level $\alpha=5\%$.
\end{itemize}
The following scenarios are considered (for the
details we refer to \citet{MaBu14}):
\begin{itemize}
\item 42 scenarios based on 7 designs;
\item for each design we consider 6 settings by varying the number of
  variables $p$ in the model and the signal to noise ratio defined by 
$\mbox{SNR} = \sqrt{(\beta^0)^T\bx^T\bx\beta^0n^{-1}\sigma^{-2}}$, namely for
$p=200$ we use $\mbox{SNR}=4$ and $\mbox{SNR}=8$, for $p=500$ we use
$\mbox{SNR}=8$ and $\mbox{SNR}=16$ and for $p=1000$ we use $\mbox{SNR}=16$ and 
$\mbox{SNR}=32$;
\item 3 designs based on synthetic data (``equi correlation'', ``high
  correlation within small blocks'' and ``high
  correlation within large blocks'') and 4 designs based on semi-real data
  (``Riboflavin with normal correlation'', ``Breast with normal
  correlation'', ``Riboflavin with high correlation'', ``Breast with high
  correlation'');
\item sparsity $s_0=6$ for the two ``Riboflavin''-designs and $s_0=10$ for the
  other five designs.
\end{itemize}


\subsection{Familywise error rate control (FWER)}

For each of the 42 scenarios described in Section \ref{sec:impscen} 
we consider exactly the same 100 independent simulation runs as in
\citet[Section 4.2.2]{MaBu14} by
varying only the synthetic noise term $\eps$ and count the number where at least
one false selection is made. According to Theorem \ref{theo:FWER}, we
expect this number to be at most $100 \alpha = 5$ ($\alpha = 0.05$). The
results for the Bonferroni-based methods can be seen in
\citet[Table 1]{MaBu14}:
FWER control holds for 40 of the 42 scenarios and in 37 scenarios there is
no false selection at all.

The results for the methods with sequential
rejection are very similar, the only differences being that for the ``high
  correlation within small
blocks''-design with $p=500$ and $\mbox{SNR}=8$ the number of runs with at
least a false selection increases (compared to Bonferroni-type methods)
from 7 to 9 for the single variable method, 
and from 7 to 13 for the hierarchical method, respectively; for the
same design 
with $p=1000$ and $\mbox{SNR}=16$ the number of runs with at least a false
selection increases from 5 to 6 for both the single variable and hierarchical
method.
For all other scenarios, inclusively the ``high correlation within large
blocks''-design with $p=200$ and $\mbox{SNR}=4$, where the non-sequential
hierarchical method slightly failed to control FWER (6 runs with at least a
false detection), the sequential methods exhibit the same FWER control
as their non-sequential counterparts.

Summarizing, FWER holds for all
four methods in 39 out of 42 scenarios and the designs where it doesn't fully
hold are ``high correlation within small blocks'' and ``high correlation
within large blocks'', 
which is not surprising since each active predictor
is highly correlated with a false variable from 
$S_0^c$ and hence it is rather difficult for our screening method (the
Lasso) to guarantee that $\hat{S} \supseteq S_0$. 

\subsection{Power}\label{sec.power}
For measuring the power we consider four different aspects: the
one-dimensional statistics defined in \citet[Section 4.2.1]{MaBu14} as
``Performance 1'' and ``Performance 2'' (see below), the number of minimal
true detections (MTDs, i.e., smallest significant groups of 
variables of any cardinality, containing at least one active variable, see below)
and singleton true detections  (STDs, i.e., MTDs with cardinality
1). Concretely, a cluster is said to be a MTD if it satisfies all of the
following:
\begin{itemize}
\item $C$ is a significant cluster, e.g., has p-value $<5\%$ (``Detection'');
\item There is no significant sub-cluster $D \subset C$ (``Minimal'');
\item $C \notin \mathcal{T}_0$, i.e., there is at least one active variable in
$C$ (``True'');
\end{itemize}
and we define: 
\begin{eqnarray*}
\mbox{Performance 1} &=& \frac{1}{|S_0|} \sum_{\mbox{MTD C}} \frac{1}{|C|},\\
\mbox{Performance 2} &=& \frac{1}{|S_0|} \sum_{\mbox{MTD C with }|C| \leq
20} \frac{1}{2} \Big( \frac{1}{|C|}+1 \Big).
\end{eqnarray*}

For each of the 42 scenarios outlined in Section \ref{sec:impscen},
we consider exactly the same 100 independent simulation runs obtained in 
\citet[Section 4.2.3-4]{MaBu14} by varying 
the synthetic noise term $\eps$ 
and the synthetic regression vector $\beta^0$. We then calculate the average 
Performance 1, Performance 2, number of MTDs and number STDs, over the 100
simulation runs. The results are shown in Table
\ref{table:powerlow} for low SNR and Table \ref{table:powerhigh} for high
SNR (for the single variable methods each MTD is an STD and by definition
Performance 2 is the same as Performance 1). 
\begin{table}[!h] \renewcommand{\tabcolsep}{4pt}
\centering
\begin{tabular}[h]{|l|c||c|c|c|c||c|c||c|c|c|c||c|c|}
\hline
& & \multicolumn{12}{c|}{low SNR}\\
\cline{3-14}
Design & $p$ & \multicolumn{4}{c||}{\# MTDs} & \multicolumn{2}{c||}{\#
  STDs} & \multicolumn{4}{c||}{Perf 1} & \multicolumn{2}{c|}{Perf 2}\\
\cline{3-14}
 & & SB & SH & HB & HSR & HB & HSR & SB & SH & HB & HSR & HB & HSR \\
\hline\hline
& 200 & 4.79 & 5.00 & 5.40 & 5.59 & 4.34 & 4.55 & 47.9 & 50.0 & 44.2 & 46.3 & 46.1 & 48.1 \\ 
\cline{2-14} 
 equi & 500 & 3.97 & 4.13 & 4.74 & 4.84 & 3.73 & 3.84 & 39.7 & 41.3 & 37.7 & 38.7 & 38.3 & 39.3 \\ 
\cline{2-14} 
 corr & 1000 & 1.77 & 1.79 & 2.54 & 2.54 & 1.73 & 1.73 & 17.7 & 17.9 & 17.4 & 17.4 & 17.6 & 17.6 \\ 
\hline 
 & 200 & 4.45 & 4.78 & 6.85 & 7.12 & 4.36 & 4.84 & 44.5 & 47.8 & 53.7 & 57.5 & 60.7 & 64.0 \\ 
\cline{2-14} 
 small & 500 & 3.15 & 3.33 & 5.18 & 5.27 & 3.15 & 3.42 & 31.5 & 33.3 & 38.3 & 40.1 & 44.2 & 45.5 \\ 
\cline{2-14} 
 blocks & 1000 & 1.31 & 1.35 & 2.53 & 2.57 & 1.31 & 1.37 & 13.1 & 13.5 & 15.1 & 15.6 & 17.2 & 17.7 \\ 
\hline 
 & 200 & 0.29 & 0.30 & 6.50 & 6.50 & 0.28 & 0.28 & 2.9 & 3.0 & 6.7 & 6.7 & 31.3 & 31.3 \\ 
\cline{2-14} 
 large & 500 & 0.06 & 0.06 & 2.76 & 2.76 & 0.06 & 0.06 & 0.6 & 0.6 & 1.1 & 1.1 & 1.1 & 1.1 \\ 
\cline{2-14} 
 blocks & 1000 & 0.00 & 0.00 & 0.60 & 0.60 & 0.00 & 0.00 & 0.0 & 0.0 & 0.1 & 0.1 & 0.1 & 0.1 \\ 
\hline 
Riboflavin & 200 & 1.41 & 1.43 & 2.41 & 2.46 & 1.33 & 1.35 & 23.5 & 23.8 & 23.4 & 23.8 & 25.2 & 25.7 \\ 
\cline{2-14} 
 normal & 500 & 0.90 & 0.90 & 1.84 & 1.85 & 0.77 & 0.79 & 15.0 & 15.0 & 13.5 & 13.8 & 14.2 & 14.5 \\ 
\cline{2-14} 
 corr & 1000 & 0.72 & 0.73 & 1.60 & 1.63 & 0.63 & 0.66 & 12.0 & 12.2 & 10.8 & 11.2 & 11.0 & 11.4 \\ 
\hline 
Breast & 200 & 4.05 & 4.16 & 5.00 & 5.11 & 3.84 & 3.94 & 40.5 & 41.6 & 39.5 & 40.6 & 41.6 & 42.9 \\ 
\cline{2-14} 
 normal & 500 & 3.95 & 4.02 & 5.04 & 5.11 & 3.82 & 3.87 & 39.5 & 40.2 & 38.8 & 39.3 & 39.6 & 40.2 \\ 
\cline{2-14} 
 corr & 1000 & 3.30 & 3.34 & 4.25 & 4.27 & 3.10 & 3.13 & 33.0 & 33.4 & 31.2 & 31.5 & 31.7 & 31.9 \\ 
\hline 
Riboflavin & 200 & 1.44 & 1.49 & 2.96 & 2.96 & 1.41 & 1.44 & 24.0 & 24.8 & 26.0 & 26.4 & 31.8 & 32.1 \\ 
\cline{2-14} 
 high & 500 & 1.72 & 1.79 & 2.95 & 2.98 & 1.69 & 1.72 & 28.7 & 29.8 & 29.9 & 30.4 & 32.8 & 33.3 \\ 
\cline{2-14} 
 corr & 1000 & 1.51 & 1.51 & 2.54 & 2.56 & 1.49 & 1.52 & 25.2 & 25.2 & 25.3 & 25.8 & 25.7 & 26.1 \\ 
\hline 
Breast & 200 & 3.98 & 4.10 & 5.91 & 5.95 & 3.87 & 3.91 & 39.8 & 41.0 & 41.2 & 41.6 & 46.1 & 46.6 \\ 
\cline{2-14} 
 high & 500 & 5.13 & 5.22 & 6.51 & 6.56 & 4.87 & 4.93 & 51.3 & 52.2 & 49.9 & 50.4 & 51.7 & 52.3 \\ 
\cline{2-14} 
 corr & 1000 & 4.73 & 4.77 & 5.95 & 5.98 & 4.64 & 4.67 & 47.3 & 47.7 & 47.0 & 47.3 & 48.3 & 48.6 \\ 
\hline 
\hline 
 \multicolumn{2}{|l||}{Average}& 2.51 & 2.58 & 4.00 & 4.06 & 2.40 & 2.48 & 27.5 & 28.3 & 28.1 & 28.8 & 31.3 & 31.9 \\ 
\hline
\end{tabular}
\caption{Number of MTDs, number of STDs, Performance 1 in \% and
  Performance 2 in \%, averaged over 100 simulation runs, for single variable
  method with Bonferroni (SB), single variable method with
Bonferroni-Holm (SH), hierarchical method with Bonferroni (HB) and
hierarchical method with sequential rejection induced by the inheritance procedure (HSR). Scenarios with low SNR.}
\label{table:powerlow}
\end{table} 
\begin{table}[!h] \renewcommand{\tabcolsep}{4pt}
\centering
\begin{tabular}[h]{|l|c||c|c|c|c||c|c||c|c|c|c||c|c|}
\hline
& & \multicolumn{12}{c|}{high SNR}\\
\cline{3-14}
Design & $p$ & \multicolumn{4}{c||}{\# MTDs} & \multicolumn{2}{c||}{\#
  STDs} & \multicolumn{4}{c||}{Perf 1} & \multicolumn{2}{c|}{Perf 2}\\
\cline{3-14}
 & & SB & SH & HB & HSR & HB & HSR & SB & SH & HB & HSR & HB & HSR \\
\hline\hline
& 200 & 9.77 & 9.83 & 9.79 & 9.80 & 9.73 & 9.74 & 97.7 & 98.3 & 97.4 & 97.5 & 97.4 & 97.5 \\ 
\cline{2-14} 
 equi & 500 & 7.28 & 7.38 & 7.63 & 7.67 & 7.18 & 7.24 & 72.8 & 73.8 & 72.0 & 72.5 & 72.1 & 72.6 \\ 
\cline{2-14} 
 corr & 1000 & 2.81 & 2.84 & 3.50 & 3.50 & 2.78 & 2.78 & 28.1 & 28.4 & 27.9 & 27.9 & 28.1 & 28.1 \\ 
\hline 
 & 200 & 9.18 & 9.31 & 9.98 & 10.00 & 9.29 & 9.48 & 91.8 & 93.1 & 96.3 & 97.4 & 98.1 & 98.7 \\ 
\cline{2-14} 
 small & 500 & 6.99 & 7.03 & 8.05 & 8.14 & 7.02 & 7.15 & 69.9 & 70.3 & 73.5 & 74.7 & 76.4 & 77.4 \\ 
\cline{2-14} 
 blocks & 1000 & 2.26 & 2.27 & 3.40 & 3.41 & 2.26 & 2.28 & 22.6 & 22.7 & 24.3 & 24.5 & 26.3 & 26.5 \\ 
\hline 
 & 200 & 2.17 & 2.26 & 9.58 & 9.58 & 2.13 & 2.14 & 21.7 & 22.6 & 27.9 & 28.0 & 61.4 & 61.4 \\ 
\cline{2-14} 
 large & 500 & 1.17 & 1.20 & 5.38 & 5.38 & 1.15 & 1.15 & 11.7 & 12.0 & 12.6 & 12.6 & 13.2 & 13.2 \\ 
\cline{2-14} 
 blocks & 1000 & 0.43 & 0.45 & 1.11 & 1.11 & 0.43 & 0.43 & 4.3 & 4.5 & 4.4 & 4.4 & 4.4 & 4.4 \\ 
\hline 
Riboflavin & 200 & 3.39 & 3.46 & 3.89 & 3.92 & 3.33 & 3.34 & 56.5 & 57.7 & 56.3 & 56.5 & 58.7 & 59.1 \\ 
\cline{2-14} 
 normal & 500 & 2.24 & 2.25 & 2.90 & 2.90 & 2.15 & 2.15 & 37.3 & 37.5 & 36.4 & 36.4 & 36.9 & 36.9 \\ 
\cline{2-14} 
 corr & 1000 & 0.98 & 1.00 & 1.83 & 1.83 & 0.96 & 0.96 & 16.3 & 16.7 & 16.2 & 16.2 & 16.3 & 16.3 \\ 
\hline 
Breast & 200 & 8.65 & 8.70 & 8.89 & 8.93 & 8.60 & 8.65 & 86.5 & 87.0 & 86.4 & 86.8 & 87.2 & 87.7 \\ 
\cline{2-14} 
 normal & 500 & 6.81 & 6.86 & 7.33 & 7.35 & 6.72 & 6.74 & 68.1 & 68.6 & 67.6 & 67.8 & 68.3 & 68.5 \\ 
\cline{2-14} 
 corr & 1000 & 3.95 & 3.97 & 4.81 & 4.84 & 3.79 & 3.82 & 39.5 & 39.7 & 38.1 & 38.4 & 38.4 & 38.7 \\ 
\hline 
Riboflavin & 200 & 3.86 & 3.97 & 4.79 & 4.83 & 3.82 & 3.86 & 64.3 & 66.2 & 66.0 & 66.7 & 69.4 & 70.1 \\ 
\cline{2-14} 
 high & 500 & 3.69 & 3.72 & 4.40 & 4.43 & 3.65 & 3.68 & 61.5 & 62.0 & 61.8 & 62.3 & 63.7 & 64.2 \\ 
\cline{2-14} 
 corr & 1000 & 2.48 & 2.51 & 3.24 & 3.27 & 2.43 & 2.45 & 41.3 & 41.8 & 40.7 & 41.0 & 41.1 & 41.4 \\ 
\hline 
Breast & 200 & 9.09 & 9.15 & 9.59 & 9.61 & 9.09 & 9.14 & 90.9 & 91.5 & 91.9 & 92.3 & 93.8 & 94.1 \\ 
\cline{2-14} 
 high & 500 & 7.75 & 7.82 & 8.38 & 8.40 & 7.71 & 7.72 & 77.5 & 78.2 & 77.8 & 78.0 & 78.9 & 79.1 \\ 
\cline{2-14} 
 corr & 1000 & 5.85 & 5.89 & 6.73 & 6.76 & 5.72 & 5.75 & 58.5 & 58.9 & 57.5 & 57.9 & 58.1 & 58.5 \\ 
\hline 
\hline 
 \multicolumn{2}{|l||}{Average}& 4.80 & 4.85 & 5.96 & 5.98 & 4.76 & 4.79 & 53.3 & 53.9 & 54.0 & 54.3 & 56.6 & 56.9 \\ 
\hline
\end{tabular}
\caption{Number of MTDs, number of STDs, Performance 1 in \% and
  Performance 2 in \%, averaged over 100 simulation runs, for single
  variable method with Bonferroni (SB), single variable method with
Bonferroni-Holm (SH), hierarchical method with Bonferroni (HB) and
hierarchical method with sequential rejection induced by the inheritance procedure (HSR). Scenarios with high SNR.}
\label{table:powerhigh}
\end{table}

Considering both low and high SNR, the methods
with sequential rejection improve the considered power measures in
comparison to the
analogous method without sequential rejection in 207 out of 252 cases, the
absolute improvement being at least 0.05 for MTDs and STDs, and at least
0.5 percent for Performance 1 and Performance 2 in 133 cases out of 
252 cases. For better interpretation of these results: an absolute
improvement of 0.05 MTDs (resp. STDs) basically means that in one out of 20
runs one more MTD (resp. STD) could be detected. Averaging over all
scenarios, the improvement given by the sequential rejection procedures
lies between 0.04 and 0.06 for MTDs and STDs, and between 0.5 and 0.7
percent for Performance 1 and Performance 2. The biggest gain with
sequential rejection can be found in the ``high correlation within small
blocks''-design with $p=200$ and low SNR: it consists of 0.48 more 
STDs, 0.33 more MTDs and an absolute increase of 3.8 percent of Performance 1 and
3.3 percent of Performance 2, respectively. This basically means that in
half of the runs the method with sequential rejection could find
one STD more and in one third of the runs it could find one MTD more.
Other particularly 
favorable scenarios for an improvement with sequential rejection are the ``equi
correlation''-design and the ``breast normal corr''-design, both with $p=200$ and
low SNR and the ``high correlation within small blocks''-design with high SNR and $p=200$, resp. $p=500$.

In general, the improvement given by the sequential rejection procedures
decreases with increasing number $p$ of covariables and is substantial only
when the power of the method without sequential rejection is intermediate. 
These empirical findings are not surprising, since looking at how
the methods are defined and in particular at the equations
(\ref{eq:originalsingle}), (\ref{eq:Holmsingle}), (\ref{eq:originalhier})
and (\ref{eq:GFhier}), we conclude that an improvement with the sequential
rejection methods is only possible if the related non-sequential method
provides at least an STD (and gets more likely the more STDs are provided by the
non-sequential method). Moreover, an improvement with sequential rejection
is more likely to happen when the number $|\hat{S}|$ of screened variables
is small.

\begin{figure}[!h]
\centerline{\includegraphics[width=1\textwidth, angle=0]{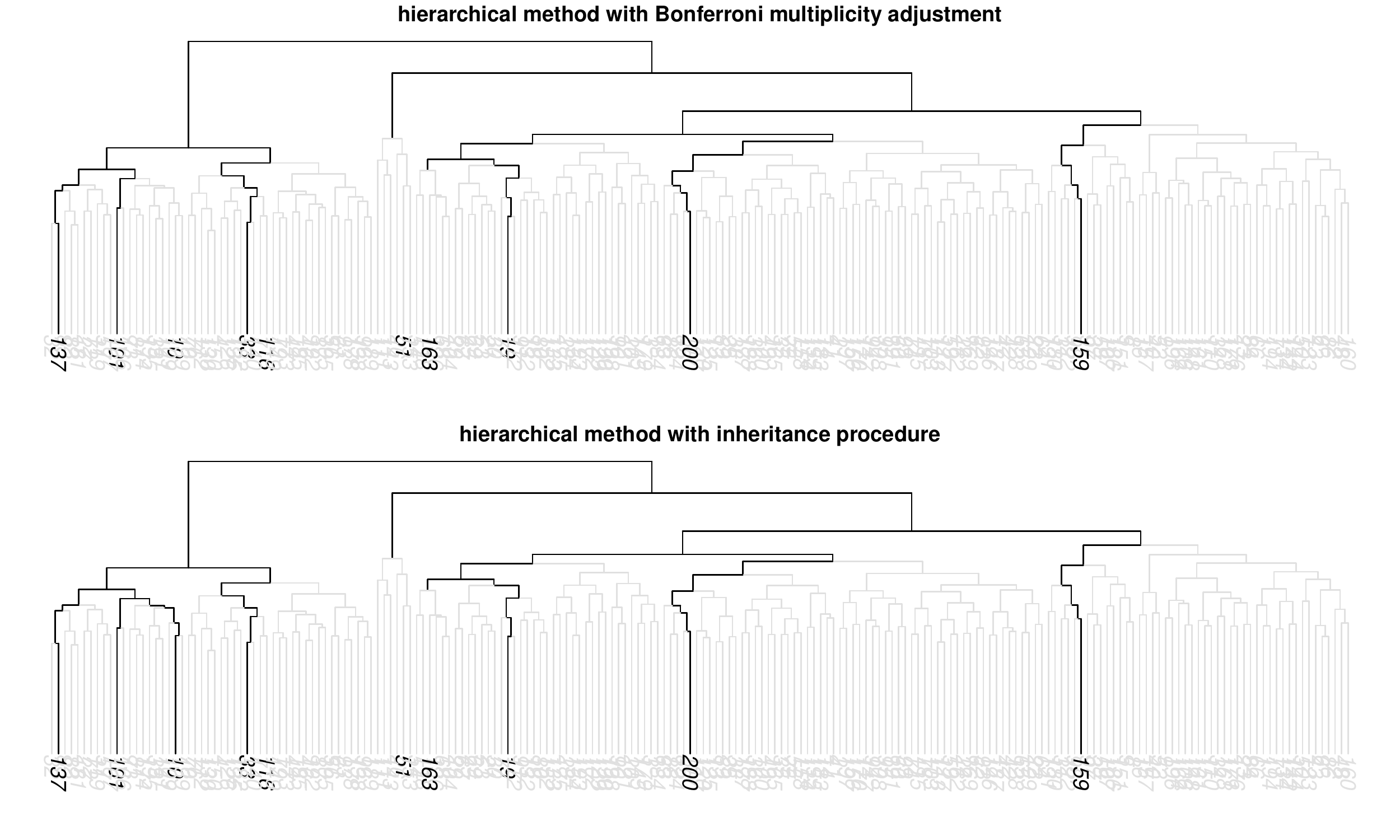}}
\vspace{-0.5cm}
\caption{Dendrograms for a paradigmatic simulation run of
  the ``equi correlation''-design with $p=200$ and $\mbox{SNR}=4$. The active
  variables are labeled in black and the truly detected 
non-zero variables along the hierarchy are depicted in black.}
\label{figure:dendro1}
\end{figure}

\begin{figure}[!h]
\centerline{\includegraphics[width=1\textwidth, angle=0]{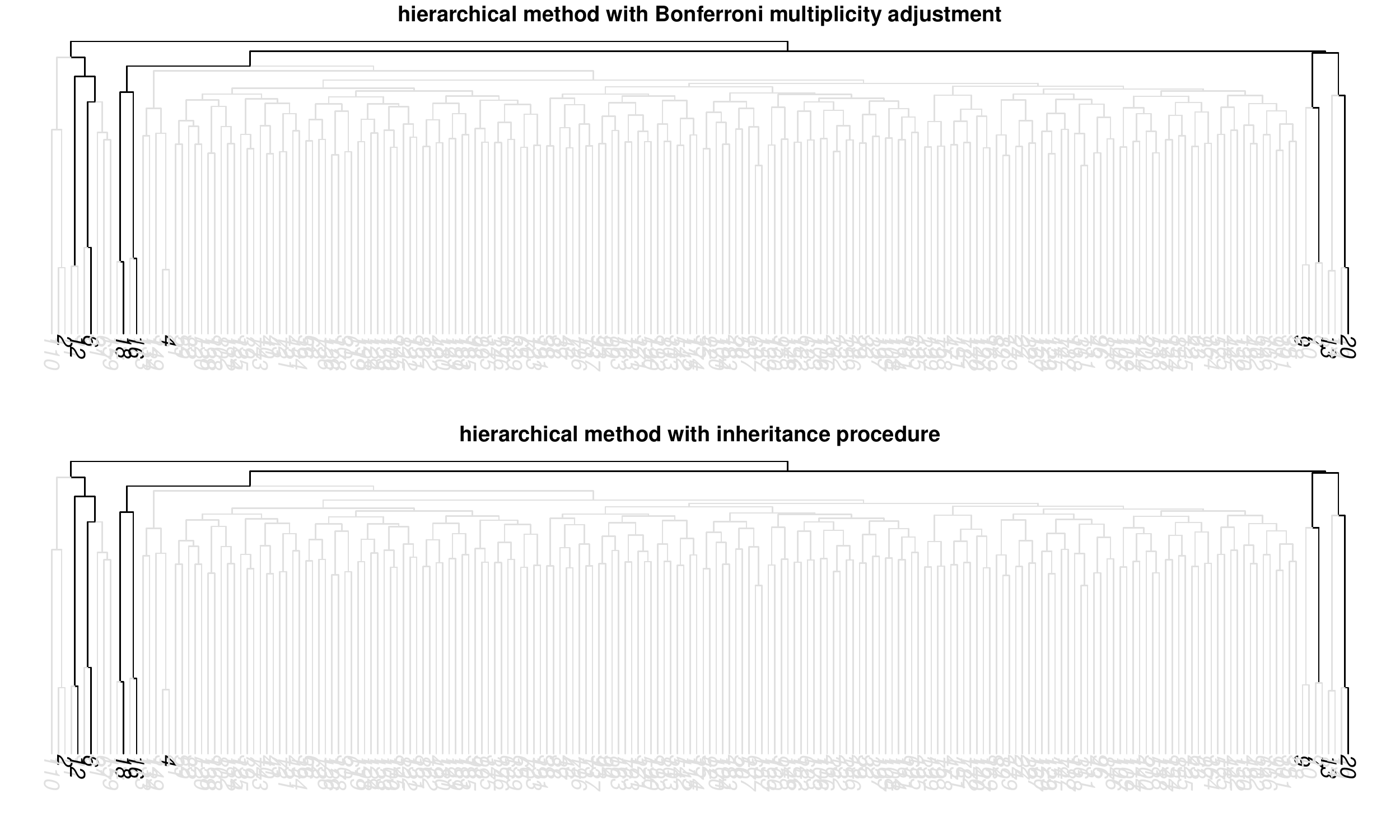}}
\vspace{-0.5cm}
\caption{Dendrograms for a paradigmatic simulation run of
  the ``high
  correlation within small blocks''-design with $p=200$ and $\mbox{SNR}=4$. The active 
  variables are labeled in black and the truly detected 
non-zero variables along the hierarchy are depicted in black.}
\label{figure:dendro2}
\end{figure}

For a better illustration of what kind of an improvement is possible using
sequential rejection, we 
show in Figures \ref{figure:dendro1} and \ref{figure:dendro2}
the dendrograms (in gray) for a paradigmatic 
simulation run of the ``equi 
  correlation''- and the ``high
  correlation within small blocks''-design, respectively, both with
  $p=200$ and $\mbox{SNR}=4$.
Figure \ref{figure:dendro1} illustrates that sequential rejection allows
the detection of a further singleton, increasing the number of STDs
from 6 to 7 and the number of MTDs 8 to 9. In Figure \ref{figure:dendro2}
sequential rejection allows to detect a singleton that could previously
only be detected together with another non-relevant variable in a cluster
of cardinality 2, increasing the number of true STDs
from 4 to 5 (while the number of MTDs remains to be 6).

Finally, we have performed a simulation with the same scenarios (and
the same sample splits) as in \citet[Section 4.3]{MaBu14}, i.e. ``small
blocks''-designs and ``large blocks''-designs with 8 different correlations
$\rho \in \{ 0, 0.4, 0.7, 0.8, 0.85, 0.9, 0.95, 0.99\}$. The full results
are shown in Tables \ref{table:corblockshighSNR} and
\ref{table:corblockslowSNR} in the Appendix. While the methods with sequential
rejection control the FWER in exactly the same scenarios where it is also
controlled by the non-sequential methods, 
they increase the average number of MTDs from 5.51
to 5.62 for the single variable method, and from 8.11 to 8.18 for the
hierarchical method, and the number of STDs for the hierarchical
method from 5.44 to 5.55, with improvements for a single scenario up to
0.48 MTDs and 0.56 STDs (averaged over 100 runs).

The empirical results can be summarized as follows. The methods with
sequential rejection essentially controls the FWER in the same way as the
non-sequential methods. Regarding power, sequential rejection allows for
improvements, to a similar extent for the 
single variable and the hierarchical procedures. As already noted in
\citet{MaBu14}, for the non-sequential methods, the hierarchical methods have similar
STDs as the single variable methods but allow for substantially more
MTDs. Thus, our proposed hierarchical method with the inheritance procedure can
be seen as the best of the considered methods.
%

\subsection{Real data application: Motif Regression}

We consider here a problem of motif regression \citep{Conlon03} from
computational biology. We apply the four methods described above, plus the
two hierarchical methods (with and without sequential rejection) using the
recently proposed canonical correlation clustering of of \citet{buru12}, to a real dataset with
$n=287$ and $p=195$, used in \citet[Section 
4.3]{Meins08} and \citet[Section 4.4]{MaBu14}. The sequential rejection
methods detects exactly the same significant structures as
non-sequential methods, namely a
single variable and a cluster containing 165 variables (the latter can be
detected only with the hierarchical method with canonical correlation
clustering). 
This can barely be considered as surprising, as with only one STD by the
non-sequential methods, further improvements by the sequential methods are
rather unlikely (see Section \ref{sec.power} for more explanation and empirical
evidence).

\section{Conclusions}\label{sec:conclusions}
We propose a general sequential rejection testing method for clusters and
single variables in a high-dimensional linear model. In presence
of high correlations among the covariables, due to serious problems of
identifiability, it is essentially mandatory to focus on detecting
significant groups of variables rather than single individual covariates. 
Our method asymptotically controls
the familywise error rate (FWER), while, as a consequence of its modular
structure, allowing for unburdened power optimization. We provide an
  implementation in the \texttt{R}-package \texttt{hdi}.

We use and study the procedure for inference of
single variables but much more importantly, for hierarchically ordered
clusters of variables. With the latter, we establish a powerful scheme for
meaningful inference in a high-dimensional regression model, much beyond
considering single variables only. 
Our presented mathematical analysis on control of the FWER and power
improvement is complemented by empirical
results based on semi-real and simulated data confirming the theoretical
results.  

\section{Appendix}\label{sec:app}
\subsection{Proof of Theorem \ref{theo:FWER}}
\begin{proof}
We show that the procedure satisfies monotonicity and single-step
conditions as required by \citet[Theorem 1]{GoSo2010}, i.e.
\begin{eqnarray}
&&\mathcal{R} \subseteq \mathcal{S} \Rightarrow \mathcal{N}(\mathcal{R})
\subseteq \mathcal{N}(\mathcal{S}) \cup \mathcal{S}\label{eq:mono}\\
&&\PP[\mathcal{N}(\mathcal{F}) \subseteq \mathcal{F}] \geq (1-\delta)^B-\alpha.\label{eq:step}
\end{eqnarray}
Assume $\mathcal{R} \subseteq \mathcal{S}$ and $C \in
\mathcal{N}(\mathcal{R})$. Then by definition
$\aggr(p^{C,(1)}m^{(1)}_C(\mathcal{R}),\dots,p^{C,(B)}m^{(B)}_C(\mathcal{R}))
\leq \alpha$. The monotonicity property (A4) of the multiplicity
adjustment and the fact that the aggregation procedure is monotone
increasing imply
$$\aggr(p^{C,(1)}m^{(1)}_C(\mathcal{S}),\dots,p^{C,(B)}m^{(B)}_C(\mathcal{S}))
\leq
\aggr(p^{C,(1)}m^{(1)}_C(\mathcal{R}),\dots,p^{C,(B)}m^{(B)}_C(\mathcal{R}))$$
and hence either $C \in \mathcal{S}$ or $C \in \mathcal{N}(\mathcal{S})$
which proves (\ref{eq:mono}).
Consider the event $$\mathcal{A}=\{\,\hat{S}^{(b)} \supseteq S_0, \forall\,b=1
\dots B \,\}$$ where all screenings are satisfied. Because of the
$\delta$-screening assumption (A2) it holds $\\P(\mathcal{A}) \geq
(1-\delta)^B$ and hence 
\begin{eqnarray*}
\PP[\mathcal{N}(\mathcal{F}) \not\subseteq
\mathcal{F}] &=& \PP[\mathcal{N}(\mathcal{F}) \not\subseteq \mathcal{F} \,|\,
\mathcal{A}] \,\PP(\mathcal{A}) + \PP[\mathcal{N}(\mathcal{F})
\not\subseteq \mathcal{F} \,|\, \mathcal{A}^c] \,\PP(\mathcal{A}^c)\\ 
&\leq& \PP[\mathcal{N}(\mathcal{F}) \not\subseteq \mathcal{F} \,|\,
\mathcal{A}]+(1-(1-\delta)^B).
\end{eqnarray*}
Since 
\begin{eqnarray*}\PP[\mathcal{N}(\mathcal{F}) \not\subseteq \mathcal{F}\,|\,
\mathcal{A}] &\leq&
\PP[\bigcup_{\mathcal{C}\setminus\mathcal{F}}\{
\aggr(p^{C,(1)}m^{(1)}_C(\mathcal{F}),\dots,p^{C,(B)}m^{(B)}_C(\mathcal{F}))\leq \alpha \}]\\
&\leq& \sum_{\mathcal{C}\setminus\mathcal{F}}
\PP[\aggr(p^{C,(1)}m^{(1)}_C(\mathcal{F}),\dots,p^{C,(B)}m^{(B)}_C(\mathcal{F}))
\leq \alpha \}]\\
&\stackrel{(A3)(A6)}{\leq}& \sum_{\mathcal{C}\setminus\mathcal{F}}
\frac{\alpha}{B} \sum_{b=1}^B \frac{1\{C \cap \hat{S}^{(b)} \neq
  0\}}{m^{(b)}_C(\mathcal{F})} = \frac{\alpha}{B} \sum_{b=1}^B \sum_{\mathcal{C}\setminus\mathcal{F}} \frac{1\{C \cap \hat{S}^{(b)} \neq 0\}}{m^{(b)}_C(\mathcal{F})}\\
&\stackrel{(A5)}{\leq}& \frac{\alpha}{B} \sum_{b=1}^B 1 \leq \alpha
\end{eqnarray*}
we conclude $\PP[\mathcal{N}(\mathcal{F}) \subseteq \mathcal{F}] = 1 -
\PP[\mathcal{N}(\mathcal{F}) \not\subseteq \mathcal{F}] \geq 1 - (\alpha +
(1-(1-\delta)^B)) = (1-\delta)^B-\alpha$ which proves (\ref{eq:step}).
\end{proof}
\subsection{Proof of Proposition \ref{theo:quantiles}}
\begin{proof}
The proof was basically given in the Appendix of \citet{memepb09}.\\
In the following we omit the function $\min \{1,\cdot\}$ from the definition of
$Q(\gamma)$ in order to simplify the notation (this is possible since the
level $\alpha$ is smaller than 1). Define for $u \in (0,1)$ the function 
$$\pi(u):= \frac{1}{B} \sum_{b=1}^B 1 \{ \tilde{p}^{(b)} \leq u \}.$$
Then it holds 
\begin{eqnarray*}
Q(\gamma) \leq \alpha &\Longleftrightarrow& q_\gamma(\tilde{p}^{(1)}/ \gamma,\dots,\tilde{p}^{(B)}/ \gamma) \leq \alpha
\Longleftrightarrow q_\gamma(\tilde{p}^{(1)},\dots,\tilde{p}^{(B)}) \leq
\alpha\gamma\\
&\Longleftrightarrow& \sum_{b=1}^B 1 \{ \tilde{p}^{(b)} \leq \alpha\gamma
\} \geq B\gamma \Longleftrightarrow \pi(\alpha\gamma) \geq \gamma.
\end{eqnarray*}
Thus,
\begin{eqnarray*}
\PP(Q(\gamma) \leq \alpha) &=& \EE (1\{ Q(\gamma) \leq \alpha \})
= \EE ( 1\{ \pi(\alpha\gamma) \geq \gamma \})
\leq \frac{1}{\gamma} \EE (\pi(\alpha\gamma))\\
&=& \frac{1}{\gamma} \EE \Big(\frac{1}{B} \sum_{b=1}^B 1 \{ \tilde{p}^{(b)} \leq
\alpha\gamma \} \Big) = \frac{1}{\gamma} \frac{1}{B} \sum_{b=1}^B \EE \Big(
1 \{ \tilde{p}^{(b)} \leq \alpha\gamma \} \Big)\\
&=& \frac{1}{\gamma} \frac{1}{B} \sum_{b=1}^B \PP(\tilde{p}^{(b)} \leq
\alpha\gamma) \leq \frac{1}{\gamma} \frac{1}{B} \sum_{b=1}^B
\frac{\alpha\gamma}{m^{(b)}}1\{C \cap \hat{S}^{(b)} \neq \emptyset\}\\ 
&=& \frac{\alpha}{B} \sum_{b=1}^B \frac{1\{C \cap \hat{S}^{(b)} \neq \emptyset\}}{m^{(b)}},
\end{eqnarray*}
where the first inequality is a consequence of the Markov inequality and
the last inequality is a consequence of the assumptions that
$P(\tilde{p}^{(b)} \leq \alpha) = P(p^{(b)}m^{(b)} \leq \alpha) \leq
\alpha/m^{(b)}$ and the definition $\tilde{p}^{(b)}=1$ for $C \cap
\hat{S}^{(b)} = \emptyset$.\\ 
For a random variable $U$ taking values in $[0,1]$,
$$\sup_{\gamma \in (\gamma_{\min},1)} \frac{1 \{U \leq
  \alpha\gamma \}}{\gamma} = \left\{ \begin{array}{ll}
0, & U \geq \alpha\\
\alpha / U, & \alpha\gamma_{\min} \leq U < \alpha\\
1 / \gamma_{\min}, & U \leq \alpha\gamma_{\min}.\\
\end{array} \right.$$
and if $U$ has an uniform distribution on $[0,1]$
\begin{eqnarray*}
\EE \Big( \sup_{\gamma \in (\gamma_{\min},1)} \frac{1 \{U \leq
  \alpha\gamma \}}{\gamma} \Big) &=& \int_0^{\alpha\gamma_{\min}}
\gamma_{\min}^{-1} dx + \int_{\alpha\gamma_{\min}}^\alpha \alpha x^{-1} dx\\
&=& \gamma_{\min}^{-1} x \big|_{x=0}^{x=\alpha\gamma_{\min}} + \alpha \log
x \big|_{x=\alpha\gamma_{\min}}^{x=\alpha}\\
&=& \alpha + \alpha (\log\alpha - \log(\alpha\gamma_{\min}))\\
&=& \alpha \big(1-\log \frac{\alpha}{\alpha\gamma_{\min}} \big)
= \alpha (1 - \log \gamma_{\min}).
\end{eqnarray*}
We apply this using as $U$ the uniform distributed 
$\tilde{p}^{(b)}/m^{(b)}=p^{(b)}$ for $C \cap S^{(b)} \neq \emptyset$ and obtain 
$$\EE \Big( \sup_{\gamma \in (\gamma_{\min},1)} \frac{1 \{\tilde{p}^{(b)}/m^{(b)} \leq \alpha\gamma \}}{\gamma} \Big) \leq \alpha (1 - \log \gamma_{\min}),$$
and similarly as above
\begin{eqnarray*}
\PP\big( \inf_{\gamma \in (\gamma_{\min},1)} Q(\gamma)
\leq \alpha \big) &=& \EE\Big( \sup_{\gamma \in (\gamma_{\min},1)} 1\{
\pi(\alpha\gamma) \geq \gamma \} \Big)\\ 
&\leq& \EE\Big( \sup_{\gamma \in (\gamma_{\min},1)} \frac{1}{B}
\sum_{b=1}^B \frac{1\{\tilde{p}^{(b)} \leq \alpha\gamma\}}{\gamma} \Big)\\
&=&\EE\Big( \sup_{\gamma \in (\gamma_{\min},1)} \frac{1}{B}
\sum_{b=1}^B \frac{ 1\{\tilde{p}^{(b)} \leq \alpha\gamma \} 1\{C \cap
  S^{(b)} \neq \emptyset \} }{\gamma} \Big)\\ 
&\leq& \frac{1}{B} \sum_{b=1}^B \EE \Big( \sup_{\gamma \in (\gamma_{\min},1)}
\frac{1 \{ \tilde{p}^{(b)}/m^{(b)} \leq \alpha\gamma/m^{(b)} \} 1\{C \cap S^{(b)} \neq \emptyset\}}{\gamma} \big)\Big)\\ 
&\leq& (1-\log \gamma_{\min}) \frac{\alpha}{B} \sum_{b=1}^B \frac{1\{C \cap S^{(b)} \neq \emptyset\}}{m^{(b)}}
\end{eqnarray*}
\end{proof}

\subsection{Additional empirical results}
\begin{table}[!h] \renewcommand{\tabcolsep}{5pt}
\centering
\begin{tabular}[h]{|c||c|c|c|c||c|c|c|c||c|c|}
\hline
$\rho$ & \multicolumn{4}{c||}{FWER} & \multicolumn{4}{c||}{\# MTDs} & \multicolumn{2}{c|}{\# STDs}\\
\cline{2-11}
 & SB & SH & HB & HSR & SB & SH & HB & HSR & HB & HSR\\
\hline\hline
\multicolumn{11}{|c|}{``small blocks''-design with high SNR}\\
0 & 0 & 0 & 0 & 0 & 9.87 & 9.89 & 9.90 & 9.90 & 9.86 & 9.86 \\ 
0.4 & 0 & 0 & 0 & 0 & 10 & 10 & 10 & 10 & 10 & 10 \\ 
0.7 & 0 & 0 & 0 & 0 & 10 & 10 & 10 & 10 & 10 & 10 \\ 
0.8 & 0 & 0 & 0 & 0 & 9.85 & 9.89 & 9.98 & 9.98 & 9.90 & 9.91 \\ 
0.85 & 0 & 0 & 0 & 0 & 9.26 & 9.38 & 9.89 & 9.92 & 9.39 & 9.53 \\ 
0.9 & 0 & 0 & 0 & 0 & 9.59 & 9.65 & 10 & 10 & 9.67 & 9.79 \\ 
0.95 & 0.21 & 0.23 & 0.21 & 0.28 & 8.36 & 8.46 & 9.82 & 9.78 & 8.36 & 8.61 \\ 
0.99 & 0.92 & 0.93 & 0.92 & 0.95 & 6.72 & 6.85 & 8.06 & 8.04 & 6.73 & 6.99 \\ 
\hline 
 \multicolumn{11}{|c|}{``large blocks''-design with high SNR} \\ 
0 & 0 & 0 & 0 & 0 & 10 & 10 & 10 & 10 & 10 & 10 \\ 
0.4 & 0 & 0 & 0 & 0 & 9.98 & 9.98 & 10 & 10 & 9.99 & 9.99 \\ 
0.7 & 0 & 0 & 0 & 0 & 5.12 & 5.35 & 9.60 & 9.60 & 5.10 & 5.12 \\ 
0.8 & 0 & 0 & 0 & 0 & 9.23 & 9.43 & 10 & 10 & 9.14 & 9.15 \\ 
0.85 & 0 & 0 & 0 & 0 & 3.86 & 4.03 & 9.98 & 9.98 & 3.84 & 3.85 \\ 
0.9 & 0 & 0 & 0 & 0 & 0.06 & 0.06 & 7.17 & 7.17 & 0.06 & 0.06 \\ 
0.95 & 0 & 0 & 0 & 0 & 1.26 & 1.29 & 9.99 & 9.99 & 1.27 & 1.28 \\ 
0.99 & 0.33 & 0.33 & 0.99 & 0.99 & 3.26 & 3.26 & 7.92 & 7.92 & 3.26 & 3.26 \\
\hline
\end{tabular}
\caption{Results of the simulation with the ``small
  blocks''- and ``large blocks''-design with high SNR (SNR=8) for 8 different
  correlations $\rho$ in the design, for single
  variable method with Bonferroni (SB), single variable method with
Bonferroni-Holm (SH), hierarchical method with Bonferroni (HB) and
hierarchical method with sequential rejection induced by the inheritance procedure (HSR).}
\label{table:corblockshighSNR}
\end{table}

\begin{table}[!h] \renewcommand{\tabcolsep}{5pt}
\centering
\begin{tabular}[h]{|c||c|c|c|c||c|c|c|c||c|c|}
\hline
$\rho$ & \multicolumn{4}{c||}{FWER} & \multicolumn{4}{c||}{\# MTDs} & \multicolumn{2}{c|}{\# STDs}\\
\cline{2-11}
 & SB & SH & HB & HSR & SB & SH & HB & HSR & HB & HSR\\
\hline\hline
\multicolumn{11}{|c|}{``small blocks''-design with high SNR}\\
0 & 0 & 0 & 0 & 0 & 9.57 & 9.69 & 9.53 & 9.63 & 9.42 & 9.53 \\ 
0.4 & 0 & 0 & 0 & 0 & 8.84 & 9.06 & 8.65 & 8.81 & 8.36 & 8.51 \\ 
0.7 & 0 & 0 & 0 & 0 & 5.87 & 6.26 & 7.28 & 7.60 & 5.65 & 6.13 \\ 
0.8 & 0 & 0 & 0 & 0 & 5.53 & 5.76 & 6.79 & 7.22 & 5.33 & 5.89 \\ 
0.85 & 0.03 & 0.04 & 0.03 & 0.04 & 2.97 & 3.08 & 5.21 & 5.56 & 2.82 & 3.14 \\ 
0.9 & 0.01 & 0.01 & 0.01 & 0.01 & 3.35 & 3.55 & 5.49 & 5.86 & 3.22 & 3.60 \\ 
0.95 & 0.46 & 0.47 & 0.46 & 0.48 & 1.02 & 1.11 & 4.04 & 4.07 & 0.9 & 0.99 \\ 
0.99 & 0.55 & 0.56 & 0.54 & 0.58 & 3.62 & 3.78 & 6.01 & 6.27 & 3.42 & 3.71 \\ 
\hline 
 \multicolumn{11}{|c|}{``large blocks''-design with low SNR} \\ 
0 & 0 & 0 & 0 & 0 & 8.42 & 8.68 & 8.38 & 8.50 & 7.98 & 8.11 \\ 
0.4 & 0 & 0 & 0 & 0 & 7.61 & 8.09 & 8.98 & 8.98 & 7.44 & 7.48 \\ 
0.7 & 0 & 0 & 0 & 0 & 0.67 & 0.71 & 5.90 & 5.91 & 0.59 & 0.59 \\ 
0.8 & 0 & 0 & 0 & 0 & 0.27 & 0.27 & 6.02 & 6.02 & 0.24 & 0.24 \\ 
0.85 & 0 & 0 & 0 & 0 & 0 & 0 & 3.38 & 3.38 & 0 & 0 \\ 
0.9 & 0 & 0 & 0.06 & 0.06 & 0.38 & 0.39 & 7.59 & 7.60 & 0.38 & 0.38 \\ 
0.95 & 0.03 & 0.03 & 0.16 & 0.16 & 0.45 & 0.45 & 8.67 & 8.68 & 0.44 & 0.44 \\ 
0.99 & 0.97 & 0.97 & 1.00 & 1.00 & 1.47 & 1.48 & 5.28 & 5.27 & 1.47 & 1.48 \\ 
\hline
\end{tabular}
\caption{Results of the simulation with the ``small
  blocks''- and ``large blocks''-design with low SNR (SNR=4) for 8 different
  correlations $\rho$ in the design, for single
  variable method with Bonferroni (SB), single variable method with
Bonferroni-Holm (SH), hierarchical method with Bonferroni (HB) and
hierarchical method with sequential rejection induced by the inheritance procedure (HSR).}
\label{table:corblockslowSNR}
\end{table}

\bibliographystyle{apalike}
\bibliography{references_4.2.2015}

\end{document}